\documentclass[11pt,twoside]{amsart}

\usepackage{amsmath}
\usepackage{amsthm}
\usepackage{amsfonts, amssymb}
\usepackage{mathrsfs}
\usepackage[all]{xy}
\usepackage{url}

\usepackage{latexsym}
\usepackage[dvips]{graphicx}

\theoremstyle{plain}

\newtheorem*{mteo}{Theorem $*$}

\theoremstyle{remark}

\theoremstyle{definition}

\begin{document}

\title[The word problem and the Aharoni-Berger-Ziv conjecture on the connectivity of independence complexes]{The word problem and the Aharoni-Berger-Ziv conjecture on the connectivity of independence complexes}

\author[J.A. Barmak]{Jonathan Ariel Barmak $^{\dagger}$}

\thanks{$^{\dagger}$ Supported by grant KAW 2005.0098 from the Knut 
and Alice Wallenberg Foundation.}

\address{Mathematics Department\\
Kungliga Tekniska h\"ogskolan\\
 Stockholm, Sweden}

\email{jbarmak@kth.se}

\begin{abstract}
For each finite simple graph $G$, Aharoni, Berger and Ziv consider a recursively defined number $\psi (G) \in \mathbb{Z}\cup \{ + \infty \}$ which gives a lower bound for the topological connectivity of the independence complex $I_G$. They conjecture that this bound is optimal for every graph. We use a result of recursion theory to give a short disproof of this claim.      
\end{abstract}

\subjclass[2000]{05C69, 55P99, 03D80, 57M05.}

\keywords{Graphs, independence complexes, topological connectivity, algorithm, computability.}

\maketitle

The map $\psi$ is defined as follows: $\psi(\emptyset)=-2$; if $G$ is a non-empty discrete graph, $\psi(G)=+ \infty$; if $G$ is non-discrete with edge set $E$, $\psi(G)=max \{ min \{\psi (G-e), \psi (G \smallsetminus e)+1 \} \ | \ e\in E\}$. Here $G-e$ denotes the subgraph of $G$ obtained by removing the edge $e$ and $G\smallsetminus e$ denotes the subgraph of $G$ induced by the vertices which are not adjacent to any of the vertices of $e$.

The \textit{independence complex} $I_G$ of a finite simple graph $G$ is the simplicial complex whose simplices are the non-empty independent subsets of vertices of $G$. From an exact sequence of \cite{Mes} (Claim 3.1) and from Van-Kampen and Hurewicz Theorems it is easy to deduce that $I_G$ is $\psi (G)$-connected \cite[Theorem 2.3]{Aha}. It is conjectured in \cite[Conjecture 2.4]{Aha} that $I_G$ is not ($\psi (G)+1$)-connected, unless it is contractible. This was proved to be true in the particular case of chordal graphs \cite{Kaw}. However we will see that the conjecture is false in general, although we will not exhibit an explicit example. The following well-known result (\cite[Corollary 3.9]{Dav}) is a consequence of the non-existence of an effective way for determining whether a group $\Gamma$ given by a finite presentation is trivial or not \cite{Adj, Rab} and a construction that associates to each presentation of $\Gamma$ a $2$-dimensional complex with fundamental group isomorphic to $\Gamma$ (see \cite{Hak} for example).

\begin{mteo}
There exists no algorithm that can decide whether a finite simplicial complex is simply connected or not.
\end{mteo}

The truthfulness of the Aharoni-Berger-Ziv Conjecture would provide an algorithm (Turing machine) capable of determining if $I_G$ is simply connected for every finite simple graph $G$ (just computing $\psi (G)$ and checking if it is positive). On the other hand, given a finite simplicial complex $K$, there is a graph $G$ such that $I_G$ is isomorphic to the first barycentric subdivision of $K$. The vertices of $G$ are the simplices of $K$ and its edges are the pairs of simplices such that none of them is a face of the other. In particular, the conjecture contradicts Theorem $*$.

\end{document}